\documentclass[11pt]{article}

\usepackage[top=1in, bottom=1in, left=1.25in, right=1.25in]{geometry}
\usepackage{enumerate}

\usepackage{amsmath}          
\usepackage{amssymb}         
\usepackage{amsthm}
\usepackage{bm}

\newtheorem{thm}{Theorem}[section]
\newtheorem{cor}[thm]{Corollary}
\newtheorem{lem}[thm]{Lemma}
\newtheorem{prop}[thm]{Proposition}
\theoremstyle{definition}
\newtheorem{defn}[thm]{Definition}

\theoremstyle{remark}
\newtheorem{rem}[thm]{Remark}
\numberwithin{equation}{section}

\begin{document}

\title{\bf An Asymptotic Formula For Counting Subset Sums Over Subgroups Of Finite Fields}
\author{Guizhen Zhu 
\\Institute For Advanced Study \\ Tsinghua University,  Beijing,  P.R. China\\zhugz08@mails.tsinghua.edu.cn
\and Daqing Wan\\Department of Mathematics \\ University of California, Irvine, CA 92697-3875, USA\\dwan@math.uci.edu}
\date{2010.12.26}
\maketitle
\begin{abstract}
Let ${\mathbb F}_q$ be the finite field of $q$ elements. Let $H\subseteq {\mathbb F}^*_q$ be a multiplicative subgroup.
For a positive integer $k$ and element $b\in {\mathbb F}_q$,  we give a sharp estimate for the number of
$k$-element subsets of $H$ which sum to $b$.
\end{abstract}


\section{Introduction}

Let ${\mathbb F }_q$ be the finite field of $q$ elements of characteristic $p$. Let $H\subseteq {\mathbb F}_q$ be a subset.
Let $1\leq k\leq |H|$ be a positive integer. For $b\in {\mathbb F}_q$,  let $M_H(k,b)$  denote the number of
$k$-element subsets $S\subseteq H$ such that
$$\sum_{ a\in S} a =b.$$
The decision version of the $k$-subset sum problem for $H$ is to determine if
$M_H(k,b)>0$. This problem arises naturally from a number of important applications in coding theory and cryptography. It is a well known
NP-complete problem,  and thus there is not much more one can say about the solution number $M_H(k,b)$ in such a generality.
The main difficulty comes from the combinatorial flexibility in choosing the subset $H$ and thus the lack of algebraic structure for the subset
$H$.  From algorithmic point of view, the dynamic algorithm \cite{C} can be used to show that the decision version of the $k$-subset sum problem 
can be solved in polynomial time if $H$ is a large subset of ${\mathbb F}_q$ in the sense that $|H|$ is of size $q^{\epsilon}$ for some constant $\epsilon>0$. 
From mathematical point of view, we are more interested in the actual value of the 
solution number $M_H(k,b)$.  Ideally, we would like to have an explicit formula or an asymptotic formula for the 
solution number $M_H(k,b)$.  This is apparently too much to hope for in general. However, 
we believe that it should be possible to obtain an asymptotic formula for the number $M_H(k,b)$ for $k$ in 
certain range if $H$ is close to a large subset
of ${\mathbb F}_q$ with certain algebraic structure. For example, it is shown in \cite{pp} that if ${\mathbb F}_q -H$ is a small set, then a good asymptotic formula
for $M_H(k, b)$ can be obtained. In addition, if $H ={\mathbb F}_q$, or ${\mathbb F}^*_q$ or any additive subgroup of ${\mathbb F}_q$,
then an explicit formula for $M_H(k, b)$ (with no error term) is obtained in \cite{pp}\cite{pq}.

When $H$ is close to a multiplicative subgroup of ${\mathbb F}_q$, the situation is more complicated as the multiplication
operation is different from the addition operation in the subset sum problem. A multiplicative subgroup is far from
the additive structure. The subset sum problem in this case becomes a highly non-linear
algebraic problem with combinatorial constraints. In this paper, we study the case that $H$ is a multiplicative subgroup of
${\mathbb  F}^*_q$ and obtain a sharp asymptotic formula for the number $M_H(k, b)$ if the index $[{\mathbb F}_q^* : H]$
is reasonably small. Our main tool is the new sieve formula from \cite{pq}  together with standard character sum arguments over finite fields.

From now on, we let $H$ be a multiplicative subgroup of ${\mathbb F}^*_q$ with index $m$. Thus, $|H| = (q-1)/m$.
The number of $k$-subsets of $H$ is ${(q-1)/m \choose k}$, and the sum could be any element $b$ of the field ${\mathbb F}_q$.
One expects that in favorable cases that the $k$-subset sums are equally distributed and thus $M_H(k, b)$ should be
about ${1\over q} {(q-1)/m \choose k}$. The key is then reduced
to estimating the error term.  Our main result is the following

\begin{thm} Let $1\leq k \leq (q-1)/m$. For $b\in {\mathbb F}^*_q$,  we have the asymptotic formula
$$\left|M_H(k,b) - {1\over q} {(q-1)/m \choose k}\right| \leq \frac{2}{\sqrt{q}}{ \sqrt{q} +k + {q\over mp} \choose k};$$
and for  $b=0$,  we have
$$\left|M_H(k,0) - {1\over q} {(q-1)/m \choose k}\right| \leq { \sqrt{q} +k + {q\over mp} \choose k},$$where $p$ is the characteristic of ${\mathbb F}_q$.
\end{thm}

Because of the obvious symmetry 
$$M_H(k, b) = M_H(|H| -k, \sum_{a\in H} a  - b),$$
we may without loss of generality assume that $k \leq |H|/2 = (q-1)/2m$.  

\begin{cor}Let $p>2$.  There is an effectively computable absolute constant $0<c<1$ such that if $m<c\sqrt{q}$ and $6\ln q<k
\leq \frac{q-1}{2m}$,  then $M_H(k,b)>0$ for all $b\in {\mathbb F}_q$.
\end{cor}

Note that even in the case that $H$ is a multiplicative subgroup of ${\mathbb F}^*_q$, we still do not know a complete polynomial time
algorithm to decide if $M_H(k,b)>0$, if $m$ is large. The result in above corollary gives a partial answer to this algorithmic question. The case that
$k$ is small (say, $k \leq  6\ln q$) can be treated using an easier Brun sieve approach \cite{CW}, but to get a non-trivial lower bound, 
one needs to assume that $m$ is significantly smaller. For example, in the case $k=2$,   one may need to assume that 
$m<q^{1/4}$ to guarantee the existence of a non-trivial ${\mathbb F}_q$-rational point on the curve $X_1^m + X_2^m = b$. 
For algorithmic purpose, the small $k$ case can often be
done by a quick exhaustive search or by using the more efficient dynamic algorithm.

To illustrate our ideas, we will also consider the following related but somewhat simpler problem of counting points on
diagonal equations with distinct coordinates. Since $H$ is a subgroup of ${\mathbb F}^*_q$ with index $m |(q-1)$,
we have $H=\{ x^m | x\in {\mathbb F}_q^*\}$.  For $0\leq k \leq q-1$, let
$N_m^*(k,b)$ denote the number of solutions of the
diagonal equation
$$x_1^m+x_2^m+\dots+x_k^m=b, $$
where the $x_i$'s are in ${\mathbb F}_q^*$ and the $x_i$'s are distinct.
There is an obvious way to compute $N_m^*(k,b)$ via the classical inclusion-exclusion principle.
Define $$X=\{(x_1,x_2,\dots,x_k)\in ({\mathbb F}_q^*)^k| x_1^m+x_2^m+\dots+x_k^m=b\}.$$
Then
\[N_m^*(k,b)=\# \{(x_1,x_2,\dots,x_k)\in X|x_i\neq x_j,i\neq j\}. \]
Let
$$X_{ij}=\{(x_1,x_2,\dots,x_k)\in X|x_i= x_j,i\neq j\}, \ \ X_{ij}^c=X - X_{ij}.$$
Applying the classical inclusion-exclusion principle, we obtain
\begin{align*}
N_m^*(k,b)&=|\bigcap_{1\leq i\leq j\leq k}X_{ij}^c|\\
&= |X|-\sum_{1\leq i\leq j\leq k}|X_{ij}|+\sum_{1\leq i\leq j\leq k\atop 1\leq s\leq t\leq k}|X_{ij}\bigcup X_{st}|-\dots+(-1)^{{k\choose 2}}|\bigcap_{1\leq i\leq j\leq k}X_{ij}|.
\end{align*}
Each term on the right side can be estimated using some basic properties of Gauss sum and Jacobi sum. The main terms are of  at most
$O(q^k)$. However, the number of terms in the above inclusion-exclusion is $2^{k\choose 2}$ which can add up to a total error term which may be greater
than the main term $O(q^k)$ as soon as  $k$ is greater than $\Omega (\sqrt{q})$.
Fortunately in \cite{pq}, J.Y.Li and D. Wan presented a new sieve for distinct coordinate counting problem which can be used for our estimation.
This sieve reduces the number of total terms from $2^{{k\choose2}}$ to $k!$, allowing us to deduce non-trivial information for $k$ as large as a fraction
of $q$ (and thus much larger than $O(\sqrt{q})$.  We will introduce their sieve briefly in Section 2. Now we state our main asymptotic formula
for the number $N_m^*(k,b)$.

\begin{thm}
For all $b\in {\mathbb F}_q^*$,  we have
\begin{equation}
\left|N_m^*(k,b)-\frac{(q-1)_k}{q}\right|\leq\frac{2}{\sqrt{q}}\left(m\sqrt{q}+k+\frac{q}{p}\right)_k,
\end{equation}
and for $b=0$, we have
$$\left|N_m^*(k,0)-\frac{(q-1)_k}{q}\right|\leq\left(m\sqrt{q}+k+\frac{q}{p}\right)_k,$$
where $(t)_k = t(t-1)\cdots (t-k+1)$ for a real number $t$.
\end{thm}

Again, we have the symmetry 
$$N_m^*(k,b) = N_m^*(q-1-k, \sum_{a\in H} a^m -b).$$
Thus, we may assume that $0\leq k \leq  (q-1)/2$.

\begin{cor} Let $p>2$. 
There is an effectively computable absolute constant $0<c<1$ such that if $m<c \sqrt{q}$ and $6\ln q<k
\leq \frac{q-1}{2}$ then $N_m^*(k,b)>0$ for all $b\in {\mathbb F}_q$.
\end{cor}

Some preliminaries will be introduced briefly in section 2. In section 3, we will illustrate the asymptotic formula for the number  $N_m^*(k,b)$. The punchline will be our asymptotic formula for $M_H(k,b)$.

\section{Preliminaries}
\subsection{Li-Wan's new sieve}

In \cite{pq},  J.Y. Li and D. Wan presented a new sieve for the distinct coordinate counting problem. We will introduce it briefly.

Let $D$ be a finite set. For a positive integer $k$, let $D^k=D\times D\times\dots\times D$ be the Cartesian product of $k$ copies of $D$. Let $X\subset D^k$. Every element $x\in X$ can be written as $x=(x_1,\dots,x_k)$ with $x_i\in D$. We are interested in counting the number of elements in $X$ with distinct coordinates, i.e., the cardinality of the set
\begin{displaymath}
\overline{X}=\{(x_1,\dots,x_k)\in X|x_i\neq x_j,i\neq j\}.
\end{displaymath}
Let $S_k$ be the symmetric group. For a given permutation $\tau\in S_k$, write its disjoint cycle product as $\tau=(i_1\cdots i_{a_1})(j_1\cdots j_{a_2})\cdots(l_1\cdots l_{a_s})$, where $a_i\geq1,1\leq i\leq s$. Define the sign of $\tau$ as $sign(\tau)=(-1)^{k-l(\tau)}$,where $l(\tau)$ is the number of disjoint cycles in $\tau$. Define
\begin{displaymath}
X_{\tau}=\{(x_1,\dots,x_k)\in X|x_{i_1}=\cdots=x_{i_{a_1}},\dots,x_{l_1}=\cdots=x_{l_{a_s}}\}.
\end{displaymath}
We have the following theorem:
\begin{thm}
\[|\overline{X}|
=\sum_{\tau\in S_k}
sign(\tau)|X_{\tau}|.\]
\end{thm}
Moreover, the group $S_k$ acts on $D^k$ by permuting its coordinates, that is
$$\tau\circ (x_1,\dots,x_k)=(x_{\tau(1)},\dots,x_{\tau(s)}).$$
If $X$ is invariant under the action of $S_k$, we call it symmetric. A permutation $\tau\in S_k$ is said to be of type $(c_1,\dots,c_k)$ if $\tau$ has exactly $c_i$ cycles of length $i$. Let $N(c_1,\dots,c_k)$ denote the number of permutations of type $(c_1,\dots,c_k)$ in $S_k$.  We have the following theorem which will be used in our main results of this paper.
 \begin{thm}\label{tt}If $X$ is symmetric, then
 \begin{displaymath}
 |\overline{X}|=\sum_{\sum ic_i=k}(-1)^{k-\sum c_i}N(c_1,\dots,c_k)|X_{\tau}|.
 \end{displaymath}
 \end{thm}
\subsection{Some combinatorial formulas}
In order to prove the main results in this paper, we need to know some combinational formulas as follows. Their proofs can be found in \cite{pq}.
Let $N(c_1,\dots,c_k)$ be the number of permutations of type $(c_1,\dots,c_k)$ in $S_k$. From \cite{ps},  we can see:
\[N(c_1,\dots,c_k)=\frac{k!}{1^{c_1}c_1!2^{c_2}c_2!\cdots k^{c_k}c_k!}.\]

\begin{lem}\label{lem1}
Define the generating function \[C_k(t_1,\dots,t_k)=\sum_{\sum ic_i=k}N(c_1,\dots,c_k)t_1^{c_1}\cdots t_k^{c_k}.\]

(1) If $t_1=\dots=t_k=q$, then
\[C_k(q,\dots,q)=(q+k-1)_k.\]

(2) If $q\geq d$, $d|(q-s)$ and $t_i=q$ for $d|i$, $t_i=s$ for $d\nmid i$, then
\begin{equation}
C_k(\overbrace{s,\dots,s}^{d-1\mbox{}},q,\overbrace{s,\dots,s}^{d-1\mbox{}},q,\dots)=k!\sum_{i=0}^{\lfloor\frac{k}{d}\rfloor}{\frac{q-s}{d}+i-1\choose \frac{q-s}{d}-1}{s+k-di-1
\choose s-1}.
\end{equation}
\end{lem}

\begin{lem}
For any giver positive integers $m,n,q$ and $l$, we have
\[\sum_{i\geq0}{l+i\choose n}{q-i\choose m}\leq{l+q+1\choose m+n+1}.\]
\end{lem}

As a corollary, we get
\begin{cor}\label{c2}
For any given positive integers $s,d,k,q$ with $q\geq s$ and $d|(q-s)$, we have

\[C_k(\overbrace{s,\dots,s}^{d-1\mbox{}},q,\overbrace{s,\dots,s}^{d-1\mbox{}},q,\dots)\leq \left(s+k+\frac{q-s}{d}-1\right)_k.\]
\end{cor}
$Proof.$ \begin{align*}
C_k(\overbrace{s,\dots,s}^{d-1\mbox{}},q,\overbrace{s,\dots,s}^{d-1\mbox{}},q,\dots)&=k!\sum_{i=0}^{\lfloor\frac{k}{d}\rfloor}{\frac{q-s}{d}+i-1\choose \frac{q-s}{d}-1}{s+k-di-1
\choose s-1}\\
&\leq k!\sum_{i=0}^{\lfloor\frac{k}{d}\rfloor}{\frac{q-s}{d}+i-1\choose \frac{q-s}{d}-1}{s+k-i-1
\choose s-1}\\
&\leq k!{s+k+\frac{q-s}{d}-1\choose k}=\left(s+k+\frac{q-s}{d}-1\right)_k.
\end{align*}
\subsection{Gauss sums and Jacobi sums}

In this subsection, we review and prove some basic properties of Gauss-Jacobi sums that are needed in our proof.
\begin{defn}
A multiplicative character on $\mathbb{F}_q^*$ is a map $\chi$ from $\mathbb{F}_q^*$ to the nonzero complex numbers $\mathbb{C}^*$
that satisfies $\chi(ab)=\chi(a)\chi(b)$ for all $a,b\in \mathbb{F}_q^*$.
We extend the definition to the whole field $\mathbb{F}_q$ by defining
\[\chi(0)=\left\{\begin{array}{ll}1,&\mbox{$\chi=1$};\\
0,&\mbox{otherwise.}\end{array}\right. \]
\end{defn}

\begin{defn}
Let $\chi$ be a multiplicative character on $\mathbb{F}_q$ and $a\in \mathbb{F}_q$.
Set
$$g_a(\chi)=\sum_{t\in {\mathbb F}_q} \chi(t)\zeta^{{\rm Tr}(at)},$$
where $\zeta=e^{2\pi i/p}$ and ${\rm Tr}$ denotes the trace map from ${\mathbb F}_q$ to ${\mathbb F}_p$.
The sum $g_a(\chi)$ is called a Gauss sum on $\mathbb{F}_q$ and we often denote $g_1(\chi)$ by $g(\chi)$.
\end{defn}
\begin{prop}\label{p1}
If $\chi\neq1$, then $|g(\chi)|=\sqrt{q}$.
\end{prop}
It's proof can be seen in \cite{pr} when $q=p$, where $p$ is a prime. For $q=p^r$, the proof is similar. 

\begin{defn}
Let $\chi_1,\dots,\chi_n$ be multiplicative characters on $\mathbb{F}_q$.  We define the following four Jacobi type sums by
\[J(\chi_1,\dots,\chi_n)=\sum_{y_1+\dots+y_n=1}\chi_1(y_1)\chi_2(y_2)\dots\chi_n(y_n).\]
\[J_0(\chi_1,\dots,\chi_n)=\sum_{y_1+\dots+y_n=0}\chi_1(y_1)\chi_2(y_2)\dots\chi_n(y_n).\]
\[J^*(\chi_1,\dots,\chi_n)=\sum_{y_1+\dots+y_n=1\atop \text{all}\  y_i\neq0}\chi_1(y_1)\chi_2(y_2)\dots\chi_n(y_n).\]
\[J^*_0(\chi_1,\dots,\chi_n)=\sum_{y_1+\dots+y_n=0\atop \text{all} \ y_i\neq0}\chi_1(y_1)\chi_2(y_2)\dots\chi_n(y_n).\]
\end{defn}
The first two sums are standard Jacobi sums. The last two sums are related to the first two sums.
These Jacobi type sums  have the following properties:
\begin{prop}\label{prop1} Let  $\chi_{i_1}=\dots=\chi_{i_e}=1$ but $\chi_{i_{e+1}}\neq1,\dots,\chi_{i_n}\neq1$. Then,

\begin{displaymath}
J(\chi_1,\chi_2,\dots,\chi_n)=\left\{\begin{array}{ll}
q^{n-1}, &\mbox{if  $e=n$,}\\
0,&\mbox{if $1\leq e<n$.}\end{array}\right.
\end{displaymath}

\begin{displaymath}
J_0(\chi_1,\chi_2,\dots,\chi_n)=\left\{\begin{array}{ll}q^{n-1}, &\mbox{if  $e=n$,}\\
0,&\mbox{if $1\leq e<n$},\\
0,&\mbox{if $e=0$, $\chi_1\chi_2\dots\chi_n\neq1,$}\\
\chi_n(-1)(q-1)J(\chi_1,\chi_2,\dots,\chi_{n-1})&\mbox{otherwise.}\end{array}\right.
\end{displaymath}

\begin{displaymath}
J^*(\chi_1,\chi_2,\dots,\chi_n)=\left\{\begin{array}{ll}
\frac{1}{q}[(q-1)^n-(-1)^n], &\mbox{if  $e=n,$}\\
(-1)^eJ(\chi_{i_{e+1}},\chi_{i_{e+2}},\dots,\chi_{i_n}),&\mbox{ if $0 \leq e<n$. }\end{array}\right.
\end{displaymath}

\begin{displaymath}
J^*_0(\chi_1,\chi_2,\dots,\chi_n)=\left\{\begin{array}{ll}
\frac{1}{q}[(q-1)^n-(q-1)(-1)^{n-1}], &\mbox{if  $e=n,$}\\
(-1)^eJ_0(\chi_{i_{e+1}},\chi_{i_{e+2}},\dots,\chi_{i_n}),&\mbox{ if  $0\leq e<n.$ }\end{array}\right.
\end{displaymath}

\end{prop}
$Proof$. If $y_1,y_2,\dots,y_{n-1}$ are chosen arbitrarily in ${\mathbb F}_q$,  then $y_n$ is uniquely determined by the condition $y_1+y_2+\dots+y_n=0$.
Thus $J_0(1,1,\dots,1)=J(1,1,\dots,1)=q^{n-1}$.
An application of  the inclusion-exclusion principle gives
\begin{align*}
 J^*(1,1,\dots,1)&=\sum_{y_1+\dots+y_n=1\atop y_i\neq0}1\\
&=q^{n-1}-{n\choose1}q^{n-2}+(-1)^2{n\choose2}q^{n-3}+\dots+(-1)^{n-1}{n\choose n-1}q^0\\
&=\frac{1}{q}[(q-1)^n-(-1)^n].
\end{align*}
Similarly,
\begin{align*}
J^*_0(1,1,\dots,1)&=\sum_{y_1+\dots+y_n=0\atop y_i\neq0}1\\
&=\sum_{s\in\mathbb{F}_q^*}\sum_{y_1+\dots+y_{n-1}=-s\atop y_i\neq0}1\\
&=(q-1)J^*(\chi_1,\dots,\chi_{n-1})\\
&=\frac{1}{q}[(q-1)^n-(q-1)(-1)^{n-1}].
\end{align*}

If $e=0$, then none of the $\chi_i$ is trivial and thus $\chi_i(0)=0$. We obtain
\[\sum_{y_1+\dots+y_n=1 }\chi_1(y_1)\chi_2(y_2)\dots\chi_n(y_n)=\sum_{y_1+\dots+y_n=1\atop y_i\neq0}\chi_1(y_1)\chi_2(y_2)\dots\chi_n(y_n),\]
as the term $\chi_1(y_1)\chi_2(y_2)\dots\chi_n(y_n)$ yields 0 if there is some $y_i=0$.  Hence,  we complete the proof of $J^*(\chi_1,\chi_2,\dots,\chi_n)=J(\chi_1,\chi_2,\dots,\chi_n)$.  The proof of $J^*_0(\chi_1,\chi_2,\dots,\chi_n)=J_0(\chi_1,\chi_2,\dots,\chi_n)$ is similar.

If $1\leq  e < n$,  without loss of generality, we may assume that $\chi_1,\chi_2,\dots,\chi_e$ are trivial and the rest are nontrivial. Then
\begin{align*}
&\sum_{y_1+\dots+y_n=0}\chi_1(y_1)\chi_2(y_2)\dots\chi_n(y_n)\\
&=q^{e-1}\sum_{y_{e+1},y_{e+2},\dots,y_{n}}\chi_{e+1}(y_{e+1})\chi_{e+2}(y_{e+2}),\dots,\chi_{n}(y_{n})\\
&=q^{e-1}\sum_{y_{e+1}}\chi_{e+1}(y_{e+1})\sum_{y_{e+2}}\chi_{e+2}(y_{e+2})\dots\sum_{y_{n}}\chi_{n}(y_{n})=0.
\end{align*}
Thus $J_0(\chi_1,\chi_2,\dots,\chi_n)=0$, and similarly  for $J(\chi_1,\chi_2,\dots,\chi_n)$. For $J^*$ Jacobi sum, we can apply the
inclusion-exclusion principle and deduce
\begin{align*}
J^*(\chi_1,\chi_2,\dots,\chi_n)&=\sum_{y_1+\dots+y_n=1\atop  y_i\neq0,1\leq i\leq n}\chi_{e+1}(y_{e+1})\dots\chi_n(y_n)\\
&=\sum_{y_1+\dots+y_n=1\atop y_i\neq0,1\leq i\leq e}\chi_{e+1}(y_{e+1})\dots\chi_n(y_n)\\
&=\sum_{y_1+\dots+y_n=1}\chi_{e+1}(y_{e+1})\dots\chi_n(y_n)-\sum_{y_2+\dots+y_n=1}\chi_{e+1}(y_{e+1})\dots\chi_n(y_n)\\
&+\dots+(-1)^e\sum_{y_{e+1}+\dots+y_n=1}\chi_{e+1}(y_{e+1})\dots\chi_n(y_n)\\
&=J(1,\dots,1,\chi_{e+1},\dots,\chi_n)+\dots+(-1)^eJ(\chi_{e+1},\dots,\chi_n)\\
&=(-1)^eJ(\chi_{e+1},\dots,\chi_n).
\end{align*}
The proof for $J_0^*(\chi_1,\chi_2,\dots,\chi_n)$ is similar.

Finally, if $e=0$, then
\[J_0(\chi_1,\chi_2,\dots,\chi_n)=\sum_s\sum_{y_1+y_2+\dots+y_{n-1}=-s}\chi_1(y_1)\chi_2(y_2)\dots\chi_{n-1}(y_{n-1})\chi_n(s)\]
We can assume $s\neq0$ in the above sum and define $y_i^{\prime}=-y_i/s$. Then
\begin{align*}
&\sum_{y_1+\dots+y_{n-1}=-s}\chi_1(y_1)\chi_2(y_2)\dots\chi_{n-1}(y_{n-1})\\
&=\chi_1\chi_2\dots\chi_{n-1}(-s)\sum_{y_1^{\prime}+y_2^{\prime}+\cdots + y_{n-1}^{\prime}=1}\chi_1(y_1^\prime)\chi_2(y_2^{\prime})\dots
\chi_{n-1}(y_{n-1}^{\prime})\\
&=\chi_1\chi_2\dots\chi_{n-1}(-s)J(\chi_1,\chi_2,\dots,\chi_{n-1}).
\end{align*}

Combining these results, we have
\begin{displaymath}
J_0(\chi_1,\chi_2,\dots,\chi_n)=\chi_1\chi_2\dots\chi_{n-1}(-1)J(\chi_1,\chi_2,\dots,\chi_{n-1})\sum_{s\neq0}\chi_1\chi_2\dots\chi_n(s).
\end{displaymath}
The last sum is 0 if $\chi_1\chi_2\dots\chi_n\neq1$ and $q-1$ if  $\chi_1\chi_2\dots\chi_n=1$. The proposition is proved.

\begin{prop}\label{p2}
If $\chi_1,\chi_2,\dots,\chi_n$ are nontrivial and $\chi_1\chi_2\dots\chi_n\neq1$, we have
\begin{displaymath}
g(\chi_1)g(\chi_2)\dots g(\chi_n)=J(\chi_1,\chi_2,\dots,\chi_n)g(\chi_1\chi_2\dots\chi_n).
\end{displaymath}
\end{prop}

\begin{cor}\label{c1}
If $\chi_1,\chi_2,\dots,\chi_n$ are nontrivial and $\chi_1\chi_2\dots\chi_n=1$, then

(1)\begin{displaymath}
g(\chi_1)g(\chi_2)\dots g(\chi_n)=\chi_n(-1)qJ(\chi_1,\chi_2,\dots,\chi_{n-1}).
\end{displaymath}

(2)\begin{displaymath}\label{e7}
J(\chi_1,\chi_2,\dots,\chi_n)=-\chi_n(-1)J(\chi_1,\chi_2,\dots,\chi_{n-1}).
\end{displaymath}
\end{cor}
$Proof.$ The proofs of Proposition and corollary above \ref{p2} can also be seen in \cite{pr} when $q=p$, and when $q=p^r$, the proofs are similar.

\begin{prop}\label{p3}
Assume that $\chi_1,\chi_2,\dots,\chi_n$ are nontrivial.

(1) If $\chi_1\chi_2\dots\chi_n\neq1$, then\[|J(\chi_1,\chi_2,
\dots,\chi_n)|=q^{(n-1)/2}.\]

(2) If $\chi_1\chi_2\dots\chi_n=1$, then
\[|J_0(\chi_1,\chi_2,
\dots,\chi_n)|=(q-1)q^{(n/2)-1}.\] and
\[|J(\chi_1,\chi_2,
\dots,\chi_n)|=q^{(n/2)-1}.\]
\end{prop}
As a corollary, if all the $\chi_i$ are nontrivial, we have
\begin{equation}\label{eq1}
|J(\chi_1,\chi_2,\dots,\chi_n)|\leq q^{(n-1)/2}.
\end{equation}
$Proof.$ Their proofs are directly from Proposition \ref{p2}, Proposition \ref{p1}, Proposition \ref{prop1} and Corollary \ref{c1}.

\section{Solutions with distinct coordinates}

Recall that  $M_H(k,b)$ denotes the number of $k$-element subsets $S\subseteq H$ such that $\sum_{a\in S}a=b$.
Namely, $M_H(k,b)$ is the number of unordered $k$-tuples $(x_1,x_2,\dots,x_k)$ with distinct $x_i\in H$ such that
\begin{equation}\label{eq2}
x_1+x_2+\dots+x_k=b.
\end{equation}
If we denote $N_H(k,b)$ be the number of ordered $k$-tuples with distinct coordinates satisfying the equation above, it's clear that $N_H(k,b)=k!M_H(k,b)$.

Note that $H=\{y^m|y\in \mathbb{F}_q^*\}$. If $(x_1,x_2,\dots,x_k)$ is a $k$-tuple satisfying the equation above, then there exist some $y_i\in \mathbb{F}_q^*$ such
that $x_i=y_i^m, 1
\leq i\leq k$, and $(y_1,y_2,\dots,y_k)$ is a $k$-tuple satisfying the following equation
\begin{equation}\label{eq3}
y_1^m+y_2^m+\dots+y_k^m=b.
\end{equation}
Let  $N^*_m(k,b)$ be the number of ordered $k$-tuples $(y_1,y_2,\dots,y_k)$ with distinct $y_i\in \mathbb{F}_q^*$ satisfying equation (\ref{eq3}).
\begin{rem}
The number of solutions with distinct coordinates in $ H$ for equation (\ref{eq2}) (i.e.$N_H(k,b)$) is not equal to the number of solutions with distinct coordinates in $\mathbb{F}_q^*$ for equation (\ref{eq3})(i.e.$N^*_m(k,b)$). However there exists a delicate relationship between them, which will be described in details later.
\end{rem}

\subsection{Estimate for $N^*_m(k,b)$ with $b\not=0$}
 As defined above,
 $$N_m^*(k,b)=\#\{(x_1,x_2,\dots,x_k)\in (\mathbb{F}^*_q)^k |x_1^m+x_2^m+\dots+x_k^m=b, \ \ x_i \text{ distinct for}\  1\leq i\leq k\}.$$
 For a positive integer $d$ and element $a\in {\mathbb F}_q$,  we shall use the following well known relation: 
\[\#\{x\in\mathbb{F}_q|x^d=a\}=\sum_{\chi^d=1}\chi(a),\]
where $\chi$ runs over all multiplicative characters of ${\mathbb F}_q$ of order dividing $d$.

In this subsection,  we will estimate $N_m^*(k,b)$ for $b\in \mathbb{F}_q^*$.

\begin{lem} Let $d_1, \cdots, d_n$ be positive integers. Define
$$N^*=\#\{(x_1,\dots,x_n)\in ({\mathbb F}_q^*)^n|a_1x_1^{d_1}+\dots+a_nx_n^{d_n}=b\},$$
where $b, a_i (1\leq i\leq n)$ are in ${\mathbb F}_q^*$. Then,  we have
\begin{equation}
\left|N^*-\frac{1}{q}[(q-1)^n-(-1)^n]\right|\leq \sum_{e=0}^{n-1}\sum_{1\leq i_{e+1}<i_{e+2}<\dots<i_n\leq n}\prod_{j=e+1}^n (d_{i_j}-1)\sqrt{q}^{n-e-1}.
\end{equation}
\end{lem}
$Proof.$ Without loss of generality, we can assume $b=1$ and $d_i|(q-1)$.
\begin{align*}
N^*&=\sum_{y_1+\cdots +y_n=1\atop y_i\neq0}\prod_{i=1}^n\#\{a_ix_i^{d_i}=y_i\}=\sum_{y_1+\cdots +y_n=1\atop y_i\neq0}\prod_{i=1}^n\sum_{\chi_i^{d_i}=1}\chi_i(\frac{y_i}{a_i})\\
&=\sum_{\chi_1^{d_1}=\dots=\chi_n^{d_n}=1}\prod_{i=1}^n\chi^{-1}(a_i)\sum_{y_1+\dots+y_n=1\atop y_i\neq0}\chi_1(y_1)\cdots\chi_n(y_n)\\
&=\frac{1}{q}[(q-1)^n-(-1)^n]
\\&+\sum_{e=0}^{n-1}
\sum_{\mbox{\tiny$\begin{array}{c}\chi_1^{d_1}=\dots=\chi_n^{d_n}=1\\
\chi_{i_1}=\dots=\chi_{i_e}=1\\
\chi_{i_{e+1}},\dots,\chi_{i_n}\neq 1\end{array}$}}\prod_{i=1}^n\chi^{-1}(a_i)J^*(\chi_{{1}},\dots,\chi_{n})\\
&=\frac{1}{q}[(q-1)^n-(-1)^n]
\\&+\sum_{e=0}^{n-1}
\sum_{\mbox{\tiny$\begin{array}{c}\chi_1^{d_1}=\dots=\chi_n^{d_n}=1\\
\chi_{i_1}=\dots=\chi_{i_e}=1\\
\chi_{i_{e+1}},\dots,\chi_{i_n}\neq 1\end{array}$}}\prod_{i=1}^n\chi^{-1}(a_i)(-1)^eJ(\chi_{i_{e+1}},\dots,\chi_{i_n}).
\end{align*}
With Proposition \ref{p3}, we have
\[\left|N^*-\frac{1}{q}[(q-1)^n-(-1)^n]\right|\leq \sum_{e=0}^{n-1}\sum_{1\leq i_{e+1}<i_{e+2}<\dots <i_n\leq n}\prod_{j=e+1}^n(d_{i_j}-1)\sqrt{q}^{n-e-1}.\]
As a corollary, if $d_i=m$, for all $1\leq i\leq k$, then
\begin{align}\label{eq4}
\left|N^*-\frac{1}{q}[(q-1)^n-(-1)^n]\right|&\leq \sum_{e=0}^{n-1}{n\choose e}(m-1)^{n-e}\sqrt{q}^{n-e-1}\\
&<  (1+(m-1)\sqrt{q})^n/\sqrt{q}.
\end{align}
\begin{thm}
\label{t3}
For all $b\in \mathbb{F}_q^*$, we have
\begin{displaymath}
\left|N_m^*(k,b)-\frac{(q-1)_k}{q}\right| < \frac{2}{\sqrt{q}}\left(m\sqrt{q}+k+\frac{q}{p}\right)_k.
\end{displaymath}
\end{thm}

$Proof$. Let  $X^*=\{(x_1,\dots,x_k)\in (\mathbb{F}_q^*)^k|x_1^m+\dots+x_k^m=b\}$.  As $X^*$ is symmetric,  we can apply
\[N_m^*(k,b)=\sum_{\sum {ic_i=k}}(-1)^{k-\sum {c_i}}N(c_1,\dots,c_k)|X_\tau^*|, \]
where $\tau$ is of type $(c_1,\dots,c_k)$, and $X_{\tau}^*=\{(x_{11},\dots,x_{kc_k)}\in (\mathbb{F}_q^*)^{\sum {c_i}}|x_{11}^m+\dots+x_{1c_1}^m+2x_{21}^m+\dots+2x_{2c_2}^m+\dots+kx_{k1}^m+\dots+kx_{kc_k}^m=b\}$.
In order to compute $N_m^*(k,b)$, we should compute $|X_{\tau}^*|$ first.
 Denote
 \[\delta_i(p)=\left\{\begin{array}{ll}
0,&p\nmid i;\\
1,&p|i
\end{array}\right.\]
and  $n=\sum c_i(1-\delta_i(p))$. Then,
\begin{align*}
|X_{\tau}^*|&= (q-1)^{\sum{c_i\delta_i(p)}}
\#\{(\cdots, x_{it_i}, \cdots )\in ({\mathbb F}_q^*)^n|
\sum\limits_{1\leq i\leq k;1\leq t_i\leq c_i,\atop  p\nmid i}i x_{it_i}^m=b\}.\\
\end{align*}
Applying (\ref{eq4}),  we have
\begin{align*}
\left||X_{\tau}^*|-\frac{1}{q}(q-1)^{\sum c_i}\right|&\leq (q-1)^{\sum{c_i\delta_i(p)}}+(q-1)^{\sum{c_i\delta_i(p)}}(1+(m-1)\sqrt{q})^{\sum c_i(1-\delta_i(p))}/\sqrt{q}\\
&<  2(q-1)^{\sum{c_i\delta_i(p)}}(1+(m-1)\sqrt{q})^{\sum c_i(1-\delta_i(p))}/\sqrt{q}.
\end{align*}
Then apply Theorem \ref{tt} and Corollary \ref{c2}, we have
\begin{align*}
\left|N_m^*(k,b)-\frac{(q-1)_k}{q}\right|& < \frac{2}{\sqrt{q}}\left((m-1)\sqrt{q}+k+\frac{q-(m-1)\sqrt{q}}{p}\right)_k
\\&\leq \frac{2}{\sqrt{q}}\left((m-1)\sqrt{q}+k+\frac{q-1}{p}\right)_k\\
&\leq \frac{2}{\sqrt{q}}\left(m\sqrt{q}+k+\frac{q}{p}\right)_k.
\end{align*}

\begin{thm}\label{tt2} Let $p>2$. 
There is an effectively computable  absolute constant $0<c<1$ such that if $m<c \sqrt{q}$ and $3\ln{4q}<k
\leq \frac{q-1}{2}$ then $N_m^*(k,b)>0$ for all $b\in \mathbb{F}_q^*$.
\end{thm}
$Proof$.  Replacing $c$ by a smaller constant if necessary,  we may assume that $m<c \sqrt{q-1}< \sqrt{q}$, 
then $(m-1)\sqrt{q}\leq m\sqrt{q-1}$. By Theorem \ref{t3},  it is sufficient to prove
\begin{displaymath}
\frac{(q-1)_k}{q}\geq \frac{2}{\sqrt{q}}\left(\frac{q-1}{p}+k+m\sqrt{q-1}\right)_k,
\end{displaymath}
that is
\begin{displaymath}
\frac{(q-1)_k}{(\frac{q-1}{p}+k+m\sqrt{q-1})_k}\geq 2\sqrt{q}.
\end{displaymath}
This holds obviously when the following inequality holds:
\begin{displaymath}
\frac{q-1}{\frac{q-1}{p}+k+m\sqrt{q-1}}\geq (4q)^{\frac{1}{2k}}.
\end{displaymath}
Since $m<c \sqrt{q-1}$ and $k\leq \frac{q-1}{2}$, it is sufficient to prove the following inequality holds:
\begin{displaymath}
\frac{q-1}{\frac{q-1}{p}+k+m\sqrt{q-1}}
\geq \frac{1}{\frac{1}{p}+\frac{1}{2}+ c}\geq (4q)^{\frac{1}{2k}}.
\end{displaymath}
Now, $p\geq 3$ and thus $1/p   + 1/2 \leq 5/6$. It is sufficient to choose a positive constant $c$ 
satisfying the inequality $c \leq \frac{1}{(4q)^{\frac{1}{2k}}}-\frac{5}{6}$.
This is possible if $(4q)^{\frac{1}{2k}}<{e}^{1/6}$, where $e$ is the natural number.  That is,  if  $k>3\ln{4q}$.   The proof is complete.

\subsection{Estimate for $N_m^*(b,0)$}
We now turn to the study of the number $N_m^*(k,b)$ when $b=0$.

\begin{lem}\label{lem4} Let $d_1, \cdots, d_n$ be positive integers. Let
$$N_0^*=\#\{(x_1,\dots,x_n)\in ({\mathbb F}_q^*)^n|a_1x_1^{d_1}+\dots+a_nx_n^{d_n}=0\},$$
where $a_i\in {\mathbb F}_q^*$.
Then,
\begin{equation}
\left|N_0^*-\frac{1}{q}(q-1)^n\right|\leq \frac{q-1}{q}+\sum_{e=0}^{n-1}\sum_{\frac{l_{i_{e+1}}}{d_{i_{e+1}}}+\dots+\frac{l_{i_{n}}}{d_{i_n}}\in \mathbb{Z}\atop
1\leq l_{i_j}\leq d_{i_j}-1}(q-1)q^{{\frac{n-e}{2}-1}}.
\end{equation}
\end{lem}
$Proof$. Without loss of generality, we can assume $d_i|(q-1)$.
\begin{align*}
N_0^*&=\sum_{y_1+\cdots +y_n=0\atop y_i\neq0}\prod_{i=1}^n\#\{a_ix_i^{d_i}=y_i\}=\sum_{y_1+\cdots +y_n=0\atop y_i\neq0}\prod_{i=1}^n\sum_{\chi_i^{d_i}=1}\chi_i(\frac{y_i}{a_i})\\
&=\sum_{\chi_1^{d_1}=\dots=\chi_n^{d_n}=1}\prod_{i=1}^n\chi^{-1}(a_i)\sum_{y_1+\dots+y_n=0\atop y_i\neq0}\chi_1(y_1)\cdots\chi_n(y_n)\\
&=\frac{1}{q}[(q-1)^n-(q-1)(-1)^{n-1}]
\\&+\sum_{e=0}^{n-1}
\sum_{\mbox{\tiny$\begin{array}{c}\chi_1^{d_1}=\dots=\chi_n^{d_n}=1\\
\chi_{i_1}=\dots=\chi_{i_e}=1\\
\chi_{i_{e+1}},\dots,\chi_{i_n}\neq 1\end{array}$}}\prod_{i=1}^n\chi^{-1}(a_i)J_0^*(\chi_{{1}},\dots,\chi_{n})\\
&=\frac{1}{q}[(q-1)^n-(q-1)(-1)^{n-1}]
\\&+\sum_{e=0}^{n-1}
\sum_{\mbox{\tiny$\begin{array}{c}\chi_1^{d_1}=\dots=\chi_n^{d_n}=1\\
\chi_{i_1}=\dots=\chi_{i_e}=1\\
\chi_{i_{e+1}},\dots,\chi_{i_n}\neq 1\\
\chi_{i_{e+1}}\chi_{i_{e+2}}\dots\chi_{i_n}=1\end{array}$}}\prod_{i=1}^n\chi^{-1}(a_i)(-1)^eJ_0(\chi_{i_{e+1}},\dots,\chi_{i_n}).
\end{align*}
By the estimation for the sums $J_0$ and $J$ in Proposition \ref{p3}, we have

$$
\left|N_0^*-\frac{1}{q}(q-1)^n\right|  \leq \frac{q-1}{q}
+\sum_{e=0}^{n-1}\sum_{\mbox{\tiny$\begin{array}{c}\frac{l_{i_{e+1}}}{d_{i_{e+1}}}+\dots+\frac{l_{i_{n}}}{d_{i_n}}\in \mathbb{Z}\\
1\leq l_{i_j}\leq d_{i_j}-1\\
1\leq i_{e+1}<\dots<i_{n}\leq n\end{array}$}}(q-1)q^{{\frac{n-e}{2}-1}}.
$$

In particular, when $d_i=m$ for all $1\leq i\leq n$, we have
\begin{align*}
\left|N_0^*-\frac{1}{q}(q-1)^n\right|&\leq\frac{q-1}{q}+\frac{q-1}{q}[\sum_{e=0}^{n-1}{n\choose e}(m-1)^{n-e}\sqrt{q}^{n-e}]\\
&=\frac{q-1}{q} (1+(m-1)\sqrt{q})^n\leq (1+(m-1)\sqrt{q})^n.
\end{align*}

\begin{thm} We have
\label{t8}
\begin{equation}
\left|N_m^*(k,0)-\frac{(q-1)_k}{q}\right|\leq\left(m\sqrt{q}+k+\frac{q}{p}\right)_k.
\end{equation}
\end{thm}

$Proof$.  Let $X_0^*=\{(x_1,\dots,x_k)\in (\mathbb{F}_q^k)^*|x_1^m+\dots+x_k^m=0\}$. As $X$ is symmetric we can apply \[N_m^*(k,0)=\sum_{\sum {ic_i=k}}(-1)^{k-\sum {c_i}}N(c_1,\dots,c_k)|X^*_{0\tau}|,\]
where $\tau$ is of type $(c_1,\dots,c_k)$, and
$X^*_{0\tau}=\{(x_{11},\dots,x_{kc_k)}\in (\mathbb{F}_q^*)^{\sum {c_i}}|x_{11}^m+\dots+x_{1c_1}^m+2x_{21}^m+\dots+2x_{2c_2}^m+\dots+kx_{k1}^m+\dots+kx_{kc_k}^m=0\}$.
In order to compute $N_m^*(k,0)$, we need to compute $X^*_{0\tau}$ first. Let $\delta_i(p), n$ be defined the same way as before. We have
\begin{displaymath}
|X^*_{0\tau}|= (q-1)^{\sum{c_i\delta_i(p)}}
\#\{(\cdots, x_{it_i}, \cdots )\in (\mathbb{F}_q^*)^n|
\sum\limits_{1\leq i\leq k;1\leq t_i\leq c_i,\atop  p\nmid i}i x_{it_i}^m=0\}.
\end{displaymath}
With the result of Lemma \ref{lem4}, we can conclude

\[\left||X^*_{0\tau}|- \frac{(q-1)^{\sum{c_i}}}{q}\right|\leq (q-1)^{\sum{c_i\delta_i(p)}}(1+(m-1)\sqrt{q})^n.\]
Applying Theorem \ref{tt} and Corollary \ref{c2}, we have
\begin{align*}
\left|N_m^*(k,0)-\frac{(q-1)_k}{q}\right|&\leq \sum_{\sum_{ic_i=k}}N(c_1,c_2,\dots,c_k)(q-1)^{\sum c_i\delta_i(p)}(1+(m-1)\sqrt{q})^n\\
&\leq \left((m-1)\sqrt{q}+k+\frac{q-1}{p}\right)_k
\\&\leq\left(m\sqrt{q}+k+\frac{q}{p}\right)_k.
\end{align*}

\begin{thm}\label{tt1}Let $p>2$. 
There is an effectively computable absolute constant $0<c<1$ such that if $m<c \sqrt{q}$ and $6\ln{q}<k
\leq \frac{q-1}{2}$ then $N_m^*(k,0)>0.$
\end{thm}
$Proof$. The proof is similar to the proof of Theorem \ref{tt2}.

\section{The subset sum problem}
Suppose $H$ is a subgroup of $\mathbb{F}_q^*$, so $H=\{x^m|x\in \mathbb{F}_q^*\}$, as $\mathbb{F}_q^*$ is a cyclic group.
Our goal is to estimate $N_H(k,b)$ ( and thus $M_H(k,b)$), which is the number of solutions with distinct coordinates in $H$ of the following equation
\begin{displaymath}
x_1+x_2+\dots+x_k=b.
\end{displaymath}
Actually,  let
$$\widetilde{X}=\{(x_1,x_2,\dots,x_k)\in H^k|x_1+x_2+\dots+x_k=b\}.$$
Then
$$N_H(k,b)=\#\{(x_1,x_2,\dots,x_k)\in \widetilde{X}|x_i \  \text{distinct}\},$$
Obviously, $\widetilde{X}$ is symmetric.
So we can apply Li-Wan's new sieve formula to compute $N_H(k,b)$.

\subsection{Estimate for $N_H(k,b)$ with $b\neq0$}

  Similar to the analysis above,  we need to compute $\widetilde{X}_{\tau}$ first,
where $\tau$ is a permutation of type $(c_1,c_2,\dots,c_k)$ in $S_k$, and
$\widetilde{X}_{\tau}=\{(x_{11},\dots,x_{1c_1},\dots,x_{kc_k})\in H^{\sum c_i}|x_{11}+\dots+x_{1c_1}+\dots+kx_{kc_k}=b\}$.  Note that $\sum ic_i=k.$
\begin{thm}\label{t6}
For any $b\in \mathbb{F}_q^*$,  we have
\begin{displaymath}
\left|N_H(k,b)-\frac{1}{q}\left(\frac{q-1}{m}\right)_k\right|\leq \frac{2}{\sqrt{q}}\left(\sqrt{q}+k+\frac{q}{mp}\right)_k.
\end{displaymath}
\end{thm}

$Proof$.  As $H=\{x^m|x\in \mathbb{F}_q^*\}$,  we can write $x_{it_i}=y_{it_i}^m$ for some $y_{it_i}\in \mathbb{F}_q^*$,
where $1\leq i\leq k, 1\leq t _i\leq c_i$.  So $\widetilde{X}_{\tau}$ equals the following
\begin{displaymath}
X^*_{\tau}=\{(y_{11},\dots,y_{1c_1},\dots,y_{kc_k})\in ({\mathbb F}_q^*)^{\sum c_i}|y^m_{11}+\dots+y^m_{1c_1}+\dots+ky^m_{kc_k}=b\}.
\end{displaymath}
Note that $y^m=(y^{\prime})^m$ iff $y=y^{\prime}\xi$, where $\xi$ is an $m$-th  root of unity.
The number of variables in $\widetilde{X}_{\tau}$ is $\sum c_i$,  so
 \[|\widetilde{X}_{\tau}|=\frac{|X^*_{\tau}|}{m^{\sum c_i}},\]
 where $|X^*_{\tau}|$
has been computed in  section 3.1.
Then $|\widetilde{X}_{\tau}|$ is given by
\begin{align*}
\frac{1}{m^{\sum c_i}}(q-1)^{\sum{c_i\delta_i(p)}}
\#\{(\cdots, x_{it_i}, \cdots )\in ({\mathbb F}_q^*)^n|
\sum\limits_{1\leq i\leq k\atop 1\leq t_i  \leq c_i,  p\nmid i}i x_{it_i}^m=0\}.
\end{align*}
Thus
\begin{displaymath}
\left||\widetilde{X}_{\tau}|-\frac{1}{q}(\frac{q-1}{m})^{\sum c_i}\right|\leq \frac{2(q-1)^{\sum c_i\delta_i(p)}(1+(m-1)\sqrt{q})^{\sum c_i(1-\delta_i(p))}}{m^{\sum c_i}\sqrt{q}}.
\end{displaymath}
As analyzed above, applying Theorem \ref{tt} and Corollary \ref{c2}, we can conclude
\begin{align*}
\left|N_H(k,b)-\frac{1}{q}\left(\frac{q-1}{m}\right)_k\right|
&\leq \frac{2}{\sqrt{q}}\sum_{\sum_{ic_i=k}}N(c_1,c_2,\dots,c_k)\left(\frac{q-1}{m}\right)^{\sum c_i\delta_i(p)}\left(\frac{1+(m-1)\sqrt{q}}{m}\right)^n\\
&\leq \frac{2}{\sqrt{q}}\left(\frac{(m-1)\sqrt{q}+1}{m}+k-1+\frac{q-1}{mp}\right)_k
\\&\leq\frac{2}{\sqrt{q}}\left(\sqrt{q}+k+\frac{q}{mp}\right)_k.
\end{align*}

As a corollary, we have
\begin{cor}\label{c4.2} For $b\in {\mathbb F}_q^*$,  we have
\begin{displaymath}
\left|M_H(k,b)-\frac{1}{q}{\frac{q-1}{m}\choose k}\right| \leq {\frac{2}{\sqrt{q}}  {\sqrt{q}+k+\frac{q}{mp}\choose k}}.
\end{displaymath}
\end{cor}
\begin{thm}\label{t4.3}Let $p>2$. 
There is an effectively computable absolute constant $0<c<1$ such that if $m<c\sqrt{q}$ and $3\ln{4q}<k
\leq \frac{q-1}{2m}$ then $M_H(k,b)>0$ for all $b\in \mathbb{F}_q^*$.
\end{thm}
$Proof$.  Replacing $c$ by a smaller constant if necessary, we may assume that  $\sqrt{q} \leq c \frac{q-1}{m}$.
 By Corollary \ref{c4.2}, it is sufficient to prove
\begin{displaymath}
\frac{1}{q}\big(\frac{q-1}{m}\big)_k\geq 2( (c+\frac{1}{2})\frac{q-1}{m}+\frac{q-1}{mp})_k/\sqrt{q}.
\end{displaymath}
That is,
\begin{displaymath}
\frac{\big(\frac{q-1}{m}\big)_k}{( (c+\frac{1}{2}) \frac{q-1}{m}+\frac{q-1}{mp})_k}\geq 2\sqrt{q}.
\end{displaymath}
This holds obviously when the following inequality holds:
\begin{displaymath}
\frac{1}{ c+ \frac{1}{2}+\frac{1}{p}}\geq (4q)^{\frac{1}{2k}},
\end{displaymath}
equivalently,
$$c \leq \frac{1}{(4q)^{\frac{1}{2k}}}-(\frac{1}{p}+ \frac{1}{2}).$$
Since $p\geq 3$, we have $1/p +1/2 \leq 5/6$. 
The existence of such a positive constant $c$ is possible if $(4q)^{\frac{1}{2k}} \leq {e}^{1/6}$,
where $e$ is the natural number. That is if  $k>3\ln{4q}$.  The proof is complete.
\subsection{Estimate for $M_H(k,0)$}
In the following,  we discuss the case when $b=0$.

Let
$$\widetilde{X_0}=\{(x_1,x_2,\dots,x_k)\in H^k|x_1+x_2+\dots+x_k=0\}.$$
Then $$N_H(k,0)=\#\{(x_1,x_2,\dots,x_k)\in \widetilde{X_0}|x_i \
\text{distinct}\}.$$
 Obviously, $\widetilde{X_0}$ is symmetric.
So we can apply the new sieve formula to compute $N_H(k,0)$.
Similar to the analysis above, we need to compute $\widetilde{X}_{0\tau}$ first,  where $\tau$ is a permutation of type $(c_1,c_2,\dots,c_k)$ in $S_k$,
and $\widetilde{X}_{0\tau}=\{(x_{11},\dots,x_{1c_1},\dots,x_{kc_k})\in H^{\sum c_i}|x_{11}+\dots+x_{1c_1}+\dots+kx_{kc_k}=0\}$. Note  that $\sum ic_i=k.$
\begin{thm}\label{t9}
\begin{displaymath}
\left|N_H(k,0)-\frac{1}{q}\left(\frac{q-1}{m}\right)_k\right|\leq\left(\sqrt{q}+k+\frac{q}{mp}\right)_k.
\end{displaymath}
\end{thm}
$Proof$.  As $H=\{x^m|x\in \mathbb{F}_q^*\}$,  we can write  $x_{it_i}=y_{it_i}^m$ for some $y_{it_i}\in \mathbb{F}_q^*$,
where $1\leq i\leq k, 1\leq t_i\leq c_i$. So  $\widetilde{X}_{0\tau}$ is related to
\begin{displaymath}
X^*_{0\tau}=\{(y_{11},\dots,y_{1c_1},\dots,y_{kc_k})\in ({\mathbb F}_q^*)^{\sum c_i}|y^m_{11}+\dots+y^m_{1c_1}+\dots+ky^m_{kc_k}=0\}
\end{displaymath}
by the formula
 \[|\widetilde{X}_{0\tau}|=\frac{|X^*_{0\tau}|}{m^{\sum c_i}},\]
where $|X^*_{0\tau}|$
has been computed in section 3.2.
It follows that $|\widetilde{X}_{0\tau}|$ is given by
\begin{align*}
\frac{1}{m^{\sum c_i}}(q-1)^{\sum{c_i\delta_i(p)}}
\#\{(\dots, x_{it_i}, \dots )\in ({\mathbb F}_q^*)^n|
\sum\limits_{1\leq i\leq k\atop 1\leq t_i \leq c_i,  p\nmid i}i x_{it_i }^m=0\}.
\end{align*}
Thus
\begin{displaymath}
\left||\widetilde{X}_{0\tau}|-\frac{1}{q}(\frac{q-1}{m})^{\sum c_i}\right|\leq \frac{(q-1)^{\sum c_i\delta_i(p)}(1+(m-1)\sqrt{q})^{\sum c_i(1-\delta_i(p))}}{m^{\sum c_i}}.
\end{displaymath}
As analyzed above, applying Theorem \ref{tt} and Corollary \ref{c2}, we can conclude
\begin{align*}
\left|N_H(k,0)-\frac{1}{q}\big(\frac{q-1}{m}\big)_k\right|
&\leq \sum_{\sum_{ic_i=k}}N(c_1,c_2,\dots,c_k)\left(\frac{q-1}{m}\right)^{\sum c_i\delta_i(p)}\left(\frac{1+(m-1)\sqrt{q}}{m}\right)^n\\
&\leq \left(\frac{(m-1)\sqrt{q}+1}{m}+k-1+\frac{q-1}{mp}\right)_k\\
&\leq \left(\sqrt{q}+k+\frac{q}{mp}\right)_k.
\end{align*}

As a corollary, we have
\begin{cor}
\begin{displaymath}
\left|M_H(k,0)-\frac{1}{q}{\frac{q-1}{m}\choose k}\right|\leq {\sqrt{q}+k-1+\frac{q}{mp}\choose k} .
\end{displaymath}
\end{cor}

\begin{thm}Let $p>2$. 
There is an effectively computable absolute constant $0<c<1$ such that if $m<c \sqrt{q}$ and $6\ln{q}<k
\leq \frac{q-1}{2m}$,  then $M_H(k,0)>0.$
\end{thm}
$Proof$.
The proof is similar to the proof of Theorem \ref{t4.3}.

\bigskip

{\bf Acknowledgments.}  This work was done while the first author was visiting  
the University of California at Irvine  with a graduate fellowship from Tsinghua 
University. She would like to thank both institutions for their hospitality and 
financial support.

\end{document}